\input amstex\documentstyle{amsppt}  
\pagewidth{12.5cm}\pageheight{19cm}\magnification\magstep1
\topmatter
\title Springer's work on unipotent classes and Weyl group representations\endtitle
\author G. Lusztig\endauthor
\address{Department of Mathematics, M.I.T., Cambridge, MA 02139}\endaddress
\thanks{Supported by NSF grant DMS-1855773.}\endthanks
\endtopmatter   
\document

\define\part{\partial}

\define\n{\notin}

\define\m{\mapsto}
\define\do{\dots}

\define\T{\times}

\redefine\i{^{-1}}

\define\ps{\psi}

\redefine\l{\lambda}

\define\kk{\bold k}

\define\CC{\bold C}

\define\EE{\bold E}

\define\RR{\bold R}

\define\ZZ{\bold Z}

\define\cb{\Cal B}

\define\fg{\frak g}

\subhead 1\endsubhead
This paper will discuss some of the contributions of T. A. Springer (1926-2011) to the theory of algebraic groups,
with emphasis on his work on unipotent classes and representations of Weyl groups. Many aspects of Springer's
work such as his contributions to the theory of Jordan algebras will not be discusses here.

Let $G$ be a connected reductive group over an algebraically closed field $\kk$ of characteristic $p\ge0$.
Let $\fg$ be the Lie algebra of $G$.
Many questions about the structure and representations of $G$ and its forms over non-algebraically closed fields
involve in essential way the understanding of unipotent elements in $G$; these elements, which fall into finitely
many conjugacy classes, have a complicated combinatorial structure. A very
important breakthrough in understanding this combinatorial structure was achieved by Springer.
This breakthrough introduced a remarkable connection between unipotent classes in $G$ and irreducible
representations of the Weyl group. This connection has become an indispensable ingredient in current work
in the representation theory of groups over $F_q,\CC,\RR$ and in the theory of character sheaves.
A $q$-analog of it gives rise to 
representations of affine Hecke algebras hence of $p$-adic groups. 

\subhead 2\endsubhead
When $G=GL_n(\CC)$, the unipotent conjugacy classes were classified by Weierstrass (1867) and Jordan (1870).
L. E. Dickson (1901,1904) classified
all conjugacy classes (in particular unipotent ones) in $Sp_4(F_q),Sp_6(F_q)$.
J. Williamson (1937, 1939) classified all conjugacy classes (in particular unipotent ones) in
$Sp_{2n}(\kk)$ with $p=0$.
Springer (Ph.D. Thesis, Leiden 1951, under H. D. Kloosterman) classifies all conjugacy classes (in particular 
unipotent ones) in $Sp_{2n}(\kk)$, with $p\ne2$; he also gives the structure of centralizers. He says:
``The classification for the other classical linear groups (orthogonal and unitary) can be done with an 
analogous method'' (this is a transation from the Dutch original).

\subhead 3\endsubhead
In \cite{St65} Steinberg proved the existence of regular unipotent elements in $G$. An independent
proof in the case where $p$ is a good prime for $G$ (semisimple) is given by Springer in \cite{S66}
where he also proves the
following beautiful result: the centralizer of a regular unipotent element equals the centre times a unipotent
group $U$ and $U$ is connected if and only if $p$ is a good prime for $G$. In this paper it is also
shown that, if $p$ is not a good prime for $G$,
a regular unipotent element is not contained in the identity component of $U$. This is based on a very delicate
analysis of the structure constants of a $\ZZ$-form of the Lie algebra of $G$. 

In \cite{S66a}, Springer
shows that the centralizer of any element of $G$ (for example a unipotent one)
contains a closed abelian  subgroup of dimension equal to the rank of $G$.

\subhead 4\endsubhead
Let $G^{un}$ (resp. $\fg^{nil}$) be the variety of unipotent (resp. nilpotent) elements in $G$ (resp. $\fg$).
Assume that $p$ is a good prime for $G$.
In \cite{S69}, Springer shows that there exists a bijective
morphism $f:G^{un}@>>>\fg^{nil}$  which commutes with the $G$-action. This is defined as follows.
Let $u$ be regular unipotent in $G$. Let $z(u)$ be the Lie algebra of the centralizer of $u$ in $G$. Pick
a regular nilpotent element $X\in z(u)$. Define a map $f_0$ from the $G$-orbit of $u$ to the
$G$-orbit of
$X$ by $gug\i\m gXg\i$. This is well defined and it extends to the required bijection $f$ (the ``Springer
bijection''). Now $f$ is not unique; it depends on a number of parameters equal to the rank of $G$. The
proof uses \cite{S66} which provides knowledge of centralizers of a regular unipotent/nilpotent
element. (Later it was shown that $f$ is an isomorphism except in type A where one needs a stronger hypothesis
on $p$.) For example, if $G=GL_n(\kk)$, then $1+e\m a_1e+a_2e^2+...+a_{n-1}e^{n-1}$
is a Springer bijection for any $(a_1,a_2,\do,a_n)\in\kk^*\T\kk\T\do\T\kk$. In 1999 Serre showed that
the map on conjugacy classes induced by a Springer bijection is canonical.

\subhead 5\endsubhead
In this subsection we assume that $\kk$ is an algebraic closure of a finite field $F_q$ and that $G$ has a
fixed $F_q$-rational structure. In \cite{G}, Green found all irreducible characters of $G(F_q)$ when $G=GL_n$.
In \cite{Sr}, Srinivasan found all irreducible characters of $G(F_q)$ when $G=Sp_4$ and $q$ is
odd.
In 1968-69 there was a seminar at the Institute for Advanced Study attended by Springer, Macdonald, Srinivasan,
Carter and others. The papers of Green and Srinvasan were studied and Macdonald formulated his conjecture
(mentioned in a paper by Springer for this seminar) which associates an irreducible character of $G(F_q)$ to a
maximal torus of $G$ defined over $F_q$ and a sufficiently general
character of it. In a 1968 manuscript Macdonald gave a conjectural
formula for its character which involved some unknown quantities. The missing part was a function on
the unipotent elements in $G(F_q)$ (or equivalently, via the Springer bijection, at least
when $p$ is a good prime, a function on the nilpotent elements in $\fg(F_q)$) which gives the character values
at unipotent elements (the ``Green functions''). Then in \cite{S70}
Springer came up with a new idea (inspired by work of Harish-Chandra): he proposed a definition of the
Green function associated to a maximal torus $T$ defined over $F_q$ (assuming that $p$ is large enough):
$$f_T(N)=(\text{\rm constant})\sum_A\ps(<N,A>)$$
where $N\in\fg(F_q)$ is nilpotent,
$A$ runs over the rational points of a $G(F_q)$-orbit of a regular element defined over $F_q$ in the Lie
algebra of $T$; $\ps$ is a non-trivial character $F_q@>>>\CC^*$, $<,>$ is the Killing form
$\fg\T\fg@>>>\kk$. At this point there was an explicit candidate for the irreducible character of
$G(F_q)$ predicted by Macdonald (with $p$ large enough). Eventually the representation iself
 was constructed (with no restriction on $p$) in \cite{DL}.
Springer's hypothesis (with $p$ large enough) was proved in \cite{K}.

\subhead 6\endsubhead
Let $W$ be the Weyl group of $G$. Let $v$ be an indeterminate and let
$H$ be the Hecke algebra over $\CC(v)$ associated to $W$. Let
$\{T_w;w\in W\}$ be the standard basis of $H$; if $s\in W$
is a simple reflection then $(T_s+v^2)(T_s-1)=0$. Let $\EE$ be a simple
$H$-module over an algebraic closure of $\CC(v)$. In \cite{BC}, Benson and
Curtis showed that $\EE$ can almost always be defined over
$\CC(v^2)$ (note that $H$ is obtained by extension of scalars from $\CC(v^2)$
to $\CC(v)$ from a $\CC(v^2)$-algebra.) Around 1972 Springer discovered an
example when this rationality property actually fails. I have heard about it from Springer during
his visit to the University of Warwick in 1973. I will explain below
Springer's argument.
Assume that $G$ is of type $E_7$ and that $\EE$ is one of the two
simple $H$-module of dimension $512$, namely the one which for $v=1$
becomes the irreducible representation $E$ of $W$ which comes 
from the Steinberg representation of $W/\pm1=Sp_6(F_2)$. 
Now the trace of a simple reflection $s_i$
of $W$ on $E$ is $0$. Hence the dimensions of the $+1,-1$
eigenspaces of $s_i$ on $E$ are $256,256$ and the action of $T_{s_i}\in H$ 
on $\EE$ has $256$ eigenvalues $v^2$ and $256$ eigenvalues $-1$, so that 
$\det(T_{s_i},\EE)=v^{512}(-1)^{256}=v^{512}$. 
Let $w_0$ be the longest element of $W$; it is a 
product of $63$ simple reflections and $T_{w_0}\in H$ is a product of $63$ 
factors of the form $T_{s_i}$. Hence
$\det(T_{w_0},\EE)=v^{63\T512}$. Now $T_{w_0}$ is central in $H$ hence it 
acts on $\EE$ as a scalar $\l$. 
We have $\det(T_{w_0},\EE)=\l^{512}$. Hence $\l^{512}=v^{63\T512}$ and 
$\l$ equals $v^{63}$ times a root of $1$. In particular, $\l\n\CC(v^2)$.
This provides the required example.

\subhead 7\endsubhead
In \cite{S76}, Springer discovered a connection between unipotent elements in $G$
and representations of $W$, with very important consequences for representation theory. There is a
precursor for this: in \cite{G}, Green remarks that 
if one takes leading coefficients in the Green polynomials of $GL_n(F_q)$
(character at a unipotent element indexed by a partition of $n$ of a series of representations corresponding
to a conjugacy class in the symmetric group $S_n$ indexed by another 
partition of $n$) one obtains the character table of $S_n$. In
particular, the number of irreducible components of top degree of a 
"Springer fibre" in $GL_n$ is the degree of an irreducible representation
of $S_n$. As I wrote in my Math. Review of Springer's paper, ``this observation [of Green]
finds a beautiful explanation in the paper under review''.

Let $\cb$ be the variety of Borel subgroups of $G$. 
For any $u\in G^{un}$ let $\cb_u=\{B\in\cb;u\in B\}$ (a
"Springer fibre"). This is a nonempty variety of dimension, say $d_u$.
Assume that $p$ is a sufficiently large prime number.
(This assumption is made so that Fourier transform on $\fg$ can be used.)
Springer shows that $W$ acts naturally on the $l$-adic
cohomology spaces $H^i(\cb_u)$. (Here $i\in\ZZ$ and $l$ is a prime $\ne p$.)
Springer shows that any irreducible representation of $W$ appears in
$H^{2d_u}(\cb_u)$ for some unipotent $u\in G$ which is in fact unique
up to conjugation. This defines a surjective map from the set of
irreducible representations of $W$ (up to conjugation) to the set of unipotent
conjugacy classes in $G$. This is essentially Springer's correspondence.
Springer shows that his Green functions (trigonometric sums) can be expressed 
in terms of the Springer representations (in all degrees).

In a later paper \cite{S78} Springer constructs his representations of $W$ in the case $p=0$;
as a corollary he deduces that the irreducible representations of $W$ can be defined over the rational
numbers (this was known earlier by a less transparent argument).

Springer conjectured that his representations of
$W$ can be defined without restriction on $p$; this was proved in \cite{L},
using intersection cohomology methods.

\subhead 8\endsubhead
Let $u\in G^{un},i\in\ZZ$. We consider the setup of no.5. In \cite{S84} Springer
shows that, when $p$ is a sufficiently large prime, $H^i(\cb_u)$
is pure of weight $i$. The proof is based on the idea that $H^i(\cb_u)$ is the cohomology of a projective
variety and at the same time it can be interpreted in terms of a Slodowy slice (which is smooth) and then using
estimates of Deligne for the eigenvalues of Frobenius. Combined with earlier results by Shoji and
Beynon-Spaltenstein, this implies that $H^i(\cb_u)=0$ for $i$ odd.  Springer conjectured that when $\kk=\CC$,
the integral homology of $\cb_u$ has no torsion and is zero in odd degrees. This conjecture was proved
in \cite{DLP}. When $p$ is a bad prime for $G$, it is not known whether
$H^i(\cb_u)=0$ for $i$ odd.


\widestnumber\key{ST65}
\Refs
\ref\key{BC}\by C.T.Benson and C.W.Curtis\paper On the degree and rationality 
of certain characters of finite Chevalley groups\jour Trans.Amer.Math.Soc.\yr1972
\vol165\pages251-273\endref
\ref\key {DLP} \by C.De Concini, G.Lusztig and C.Procesi\paper Homology of the zero set of a nilpotent 
vector field on a flag manifold \jour J. Amer. Math. Soc.\vol1\yr1988\pages15-34\endref
\ref\key {DL} \by P.Deligne and G.Lusztig\paper Representations of reductive groups over finite fields\jour
 Ann. Math.\vol103\yr1976\pages 103-161\endref
 \ref\key{G}\by J.A.Green\paper The characters of the finite general linear groups
 \jour Trans.Amer.Math.Soc.\yr1955\pages402-447\endref
 \ref\key{K}\by D.Kazhdan\paper Proof of Springer's hypothesis\jour Isr.J.Math.\yr1977\vol28
 \pages272-286\endref
\ref\key {L} \by G.Lusztig\paper
Green polynomials and singularities of unipotent classes \jour Adv. Math.\vol42
\yr1981\pages 169-178\endref
\ref\key{S66}\by T.A.Springer\paper Some arithmetical results on semisimple Lie algebras
\jour Inst.Hautes \'Etudes Sci.Publ.Math.\vol30\yr1966\pages 115-141\endref
\ref\key{S66a}\by T.A.Springer\paper A note on centralizers in semisimple groups\jour
Indag.Math.\vol28\yr1966\pages75-77\endref
\ref\key{S69}\by T.A.Springer\paper The unipotent variety
of a semisimple group\inbook Algebraic Geometry \bookinfo Internat.Colloq., Tata Inst.Fund.Res. Bombay
\publ Oxford Univ.Press \yr1969\pages373-391\endref
\ref\key{S71}\by T.A.Springer\paper Generalization of Green polynomials\inbook
Representations of finite groups and related topics\bookinfo Proc.Symp.Pure Math.\vol XXI
Amer.Math.Soc.\yr1971\pages149-153\endref
\ref\key{S76}\by T.A.Springer\paper Trigonometric sums, Green functions of
finite groups and representations of Weyl groups\jour Invent.Math.\yr1976
\vol36 \pages173-207\endref
\ref\key{S78}\by T.A.Springer\paper A construction of representations
of Weyl groups\jour Invent.Math.\yr1978 \vol44\pages279-294\endref
\ref\key{S80}\by T.A.Springer\paper The Steinberg function
of a finite Lie algebra\jour Invent.Math.\vol58\yr1980\pages211-215\endref
\ref\key{S84}\by T.A.Springer\paper A purity result for fixed point varieties
in flag manifolds\jour J. Fac. Sci. Univ. Tokyo IA Math.\vol31\yr1984\pages271-282\endref
\ref\key{Sr}\by B.Srinivasan\paper The characters of the finite symplectic group $Sp(4,q)$
\jour Trans. Amer. Math. Soc.\yr1968\vol131\pages488-525\endref
\ref\key{St51}\by R.Steinberg\paper A geometric approach to the representations of the
full linear group over a Galois field\jour Trans.Amer.Math.Soc.\yr1951
\vol71\pages275-282\endref
\ref\key{St65}\by R.Steinberg\paper Regular elements in semisimple algebraic groups
\jour Inst. Hautes \'Etudes Sci. Publ. Math.\vol25\yr1965\pages49-80\endref
\endRefs
\enddocument